\documentclass[11pt,reqno]{amsart}
\usepackage{graphicx}
\usepackage{verbatim}
\usepackage{textcomp}
\usepackage{amssymb}
\usepackage{cite}
\usepackage{amsmath}
\usepackage{latexsym}
\usepackage{amscd}
\usepackage{amsthm}
\usepackage{mathrsfs}
\usepackage{bm}
\usepackage{url}
\usepackage{hyperref}
\usepackage{bookmark}
\allowdisplaybreaks[3]
\newtheorem{thm}{Theorem}[section]
\newtheorem{corr}[thm]{Corollary}

\newtheorem{prop}[thm]{Proposition}

\theoremstyle{definition}

\newtheorem*{ack}{Acknowledgment}
\theoremstyle{remark}
\newtheorem{rem}{Remark}[section]
\numberwithin{equation}{section}
\begin{document}
\title[A Reilly type integral formula]
{A Reilly type integral formula and its applications}

\author{Guangyue Huang }

\address{Department of Mathematics,
Henan Normal University, Xinxiang 453007, P.R. China}
\email{hgy@htu.edu.cn(G. Huang) }

\author{Bingqing Ma }
\email{bqma@htu.edu.cn(B. Ma) }

\author{Mingfang Zhu }
\email{mfzhu21@126.com(M. Zhu) }

\thanks{The research of authors is supported by NSFC(No. 11971153).}

\begin{abstract}
In this paper, we achieve a Reilly type integral formula associated with the $\phi$-Laplacian. As its applications, we obtain
Heintze-Karcher and Minkowski type inequalities. Furthermore, almost Schur lemmas are also given. They recover the partial results of Li and Xia in \cite{lixia2019}. On the other hand, we also study eigenvalue problem for Wentzell boundary conditions and obtain eigenvalue relationships.
\end{abstract}



\maketitle

\section{Introduction}
Let $(M,g)$ be an $n$-dimensional Riemannian manifold with the dimension $n\geq3$, where $g$ is the metric.
The $\phi$-Laplacian associated with $\phi$ is defined by
$$\Delta_\phi v=e^\phi{\rm div}(e^{-\phi}\nabla v)=\Delta v-\nabla\phi\nabla v,\ \ \ \ \ \forall \ v\in
C^\infty(M),$$
which is symmetric with respect to the $L^2(M)$ inner product under the weighted measure
$$d\mu=e^{-\phi}dv_g,$$
that is,
$$\int_{M}\alpha\Delta_\phi \beta\,d\mu=\int_{M}\beta\Delta_\phi \alpha\,d\mu=-\int_{M}\nabla\alpha\nabla\beta\,d\mu, \ \ \ \forall \ \alpha,\beta\in
C^\infty_0(M).$$
Following \cite{Bakry85,li05,Wei09}(or see \cite{LW2015,HL2013,HM2016,HZ2016} and the references therein), the $m$-dimensional Bakry-\'{E}mery Ricci curvature associated
with the above $\phi$-Laplacian is given by
\begin{align}\label{Int-1}
{\rm Ric}_{\phi,m}={\rm Ric}+\nabla^2\phi-\frac{1}{m-n}d\phi\otimes d\phi,
\end{align}
where $m$ is a real constant, and $m=n$ if and only if $\phi$ is a constant. Define
$${\rm Ric}_{\phi}={\rm Ric}+\nabla^2\phi.$$
Then ${\rm Ric}_{\phi}$ can be seen as
the $\infty$-dimensional Bakry-\'{E}mery Ricci curvature, that is, ${\rm Ric}_{\phi}:={\rm Ric}_{\phi,\infty}$.

For the convenience, we still denote $\Delta,\nabla$ by the Laplacian operator and gradient operator on
$M$, and $\overline{\Delta},\overline{\nabla}$, respectively, by the Laplacian operator and gradient operator on the boundary $\partial M$. The mean curvature $H$ of $\partial M$ is given by
$H={\rm tr}_{g}(II)$, where $II(X,Y)=g(\nabla_X\nu,Y)$ denotes the second fundamental form of $\partial M$
with $\nu$ the outward unit normal on $\partial M$. For any {\em positive} twice differentiable function $V$, we make the following conventions:
\begin{equation}\label{1-Int-1}
{\rm \widehat{Ric}}^{V}_{\phi,m}=\frac{\Delta_\phi V}{V}g_{ij}-\frac{1}{V}V_{ij}+{\rm Ric}_{\phi,m},
\end{equation}
\begin{equation}\label{1-Int-2}
{\rm \widehat{Ric}}^{V}_{\phi,\infty}=\frac{\Delta_\phi V}{V}g_{ij}-\frac{1}{V}V_{ij}+{\rm Ric}_{\phi},
\end{equation}
$H_\phi=H-\phi_\nu$, $II^V=II-(\ln V)_\nu \overline{g}$ and $d\sigma$ denotes the measure induced on $\partial M$.

First, we prove the Reilly type formula associated with $\phi$-Laplacian:

\begin{thm}\label{1-thm1}
Let $V$ be a positive twice differentiable function on a given compact Riemannian manifold $M$ with the boundary $\partial M$.
For any smooth function $f$, we have the following equality:
\begin{align}\label{1-th-111}
0=&\int_{\partial M}\Big[-VII^V\Big(V\overline{\nabla}\frac{f}{V},V\overline{\nabla}\frac{f}{V}\Big)-V^3H_\phi\Big(\Big(\frac{f}{V}\Big)_\nu\Big)^2\notag\\
&-2V^2\Big(\frac{f}{V}\Big)_\nu\Big(\overline{\Delta}_\phi f-\frac{\overline{\Delta}_\phi V}{V} f\Big)\Big]\,d\sigma+\int_M V\Big[\Big(\Delta_\phi f-\frac{\Delta_\phi V}{V}f\Big)^2\notag\\
&-\Big|\nabla^2f-\frac{f}{V} \nabla^2V\Big|^2-{\rm \widehat{Ric}}^{V}_{\phi,\infty}\Big(V\nabla\frac{f}{V},V\nabla\frac{f}{V}\Big)\Big]\,d\mu.
\end{align}
\end{thm}

Denote by $A_{ij}=f_{ij}-\frac{f}{V}V_{ij}$ and $\mathring{A}_{ij}=A_{ij}-\frac{{\rm tr}(A_{ij})}{n}g_{ij}$ with ${\rm tr}(A_{ij})=\Delta f-\frac{f}{V}\Delta V$. Then, we have ${\rm tr}(\mathring{A}_{ij})=0$ and
$$\aligned
|A_{ij}|^2=&|\mathring{A}_{ij}|^2+\frac{1}{n}[{\rm tr}(A_{ij})]^2\\
=&|\mathring{A}_{ij}|^2+\frac{1}{m}[{\rm tr}(A_{ij})]^2+\frac{m-n}{mn}[{\rm tr}(A_{ij})]^2\\
=&|\mathring{A}_{ij}|^2+\frac{1}{m}\Big(\Delta_\phi f-\frac{f}{V}\Delta_\phi V\Big)^2+\frac{2}{m}V[{\rm tr}(A_{ij})]\nabla \phi\nabla \frac{f}{V} \\
&-\frac{1}{m}V^2\Big(\nabla \phi\nabla \frac{f}{V}\Big)^2+\frac{m-n}{mn}[{\rm tr}(A_{ij})]^2,
\endaligned$$
which gives
\begin{align}\label{Int-3}
&\Big|\nabla^2f-\frac{f}{V} \nabla^2V\Big|^2+{\rm Ric}_{\phi}\Big(V\nabla\frac{f}{V}, V\nabla\frac{f}{V}\Big)\notag\\
=&|\mathring{A}_{ij}|^2+\frac{1}{m}\Big(\Delta_\phi f-\frac{f}{V}\Delta_\phi V\Big)^2+{\rm Ric}_{\phi,m}\Big(V\nabla\frac{f}{V}, V\nabla\frac{f}{V}\Big)\notag\\
&+\Big(\sqrt{\frac{m-n}{mn}}{\rm tr}(A_{ij})+\sqrt{\frac{n}{m(m-n)}}V\nabla \phi\nabla \frac{f}{V}\Big)^2\notag\\
\geq&\frac{1}{m}\Big(\Delta_\phi f-\frac{f}{V}\Delta_\phi V\Big)^2+{\rm Ric}_{\phi,m}\Big(V\nabla\frac{f}{V}, V\nabla\frac{f}{V}\Big)
\end{align}
provided $m\in (-\infty,0)\cup [n,+\infty)$, and equality holds if and only if
\begin{align}\label{Int-4}
\nabla^2 f-\frac{f}{V}\nabla^2 V=\frac{1}{n}\Big(\Delta f-\frac{f}{V}\Delta V\Big)g
\end{align}
and
\begin{align}\label{Int-5}
\Delta f-\frac{f}{V}\Delta V+\frac{n}{\sqrt{(m-n)^2}}V\nabla \phi\nabla \frac{f}{V}=0.
\end{align}
Applying \eqref{Int-3} in \eqref{1-th-111}, one obtain the following result immediately:

\begin{corr}\label{1-corr2}
Let $V$ and $f$ be as in Theorem \ref{1-thm1} and $m\in (-\infty,0)\cup [n,+\infty)$. Then, we have the following inequalities:
\begin{align}\label{2-corr-111}
0\leq&\int_{\partial M}\Big[-VII^V\Big(V\overline{\nabla}\frac{f}{V},V\overline{\nabla}\frac{f}{V}\Big)-V^3H_\phi\Big(\Big(\frac{f}{V}\Big)_\nu\Big)^2\notag\\
&-2V^2\Big(\frac{f}{V}\Big)_\nu\Big(\overline{\Delta}_\phi f-\frac{\overline{\Delta}_\phi V}{V} f\Big)\Big]\,d\sigma+\int_M V\Big[\frac{m-1}{m}\Big(\Delta_\phi f-\frac{\Delta_\phi V}{V}f\Big)^2\notag\\
&-{\rm \widehat{Ric}}^{V}_{\phi,m}\Big(V\nabla\frac{f}{V},V\nabla\frac{f}{V}\Big)\Big]\,d\mu,
\end{align}
where the equality occurs if and only if \eqref{Int-4} and \eqref{Int-5} hold.

\end{corr}

\begin{rem}
In \cite{QiuXia2015} (or see \cite{Xia2016,lixia2019}), Qiu and Xia provide a generalization of Reilly's formula and give some applications. In particular, if $m=n$, then our \eqref{1-th-111} becomes the formula (1.1) of Li and Xia in \cite{lixia2019}.

\end{rem}

Next, by using the above formula \eqref{2-corr-111}, we can achieve the following Heintze-Karcher type inequality and Minkowski type inequality:

\begin{thm}\label{1-thm3}
Let $V$ be a positive twice differentiable function on a given compact Riemannian manifold $M$ with the boundary $\partial M$. If ${\rm \widehat{Ric}}^{V}_{\phi,m}\geq0$ and $H_\phi>0$, where $m\in (-\infty,0)\cup [n,+\infty)$, then
\begin{align}\label{3-th-111}
\int_MV\,d\mu\leq\frac{m-1}{m}\int_{\partial M}\frac{V}{H_\phi}\,d\sigma.
\end{align}
Moreover, if the equality in \eqref{3-th-111} holds, then $m=n$ and $\partial M$ is umbilical.
\end{thm}

\begin{thm}\label{1-thm4}
Let $V$ be a positive twice differentiable function on a given compact Riemannian manifold $M$ with the boundary $\partial M$. If ${\rm \widehat{Ric}}^{V}_{\phi,m}\geq0$ and
\begin{align}\label{4-th-111}
II^V\geq0,
\end{align}
where $m\in (-\infty,0)\cup [n,+\infty)$, then
\begin{align}\label{4-th-222}
\Big(\int_{\partial M}V\,d\sigma\Big)^2\geq\frac{m}{m-1}\int_MV\,d\mu\int_{\partial M} V H_\phi \,d\sigma.
\end{align}
Moreover, if the inequality in \eqref{4-th-111} is strict and the equality in \eqref{4-th-222} holds, then we have that $m=n$, $\partial M$ is umbilical and  the mean curvature $H$ is constant.
\end{thm}

\begin{rem}
If $\partial M=\Sigma\cup(\cup_{l=1}^\tau N_l)$ satisfies the inner boundary condition defined by Definition 2.3 in \cite{lixia2019} and we suppose the outermost boundary hypersurface $\Sigma$ is mean convex, the similar results as in Theorem \ref{1-thm3} can also be achieved. Particularly, when $m=n$, then our Theorem \ref{1-thm3} and Theorem \ref{1-thm4} reduce to the Theorem 1.3 and Theorem 1.7 of Li and Xia in \cite{lixia2019}, respectively.

\end{rem}

Let $\mathbb{R}^{n+1}(c)$ be a space form with constant sectional curvature $c$, where $\mathbb{R}^{n+1}(1)=\mathbb{S}^{n+1}(1)$ is a unit sphere if $c=1$, $\mathbb{R}^{n+1}(-1)=\mathbb{H}^{n+1}(-1)$ is a hyperbolic space if $c=-1$ and $\mathbb{R}^{n+1}(0)=\mathbb{R}^{n+1}$ is an Euclidean space if $c=0$. Let
\begin{align}\label{const-1}
C_{n,K_1,K_2,\eta_1}=\frac{n-1}{n}+\frac{(n-1)K_1+K_2+2\sqrt{[\eta_1+(n-1)K_1]K_2}}{\lambda_1},
\end{align}
where $\eta_1$ is the first nonzero eigenvalue of problem $\Delta_\phi f-\frac{\Delta_\phi V}{V}f=-\eta \frac{f}{V}$ on closed manifolds. Next, we will apply Theorem \ref{1-thm1} to obtain the following almost Schur lemmas on closed hypersurfaces:

\begin{thm}\label{1-thm5}
Let $M$ be an $n$-dimensional closed hypersurface in a space form $\mathbb{R}^{n+1}(c)$ and $V$ be a positive twice differentiable function on $M$. If $V{\rm \widehat{Ric}}^{V}_{\phi,\infty}\geq-(n-1)K_1$, where $K_1$ is a nonnegative constant, then
\begin{align}\label{5-th-111}
\frac{(n-1)^2}{n^2}\int_M V(H-\overline{H}^V)^2\,d\mu
\leq C_{n,K_1,K_2,\eta_1}\int_M V\Big(|II|^2-\frac{1}{n}H^2\Big)\,d\mu,
\end{align}
where $\overline{H}^V=\frac{\int_M VH\,d\mu}{\int_M V\,d\mu}$ and $K_2=\max(V|\nabla \phi|^2)$ for some nonnegative constant $K_2$.

\end{thm}

The $r$-mean curvature $S_r$ on hypersurface $M$ of $\mathbb{R}^{n+1}(c)$ is related to the Newton transformation $P_r$ by
$${\rm tr}(P_r)=(n-r)S_r,$$
where
$$
P_{rij}=\frac{1}{r!}\sum_{\mbox{\tiny$\begin{array}{c}
i_1,\cdots,i_{r} \\
j_1,\cdots,j_{r}\end{array}$}}\delta_{i_1\cdots i_r i}^{j_1\cdots
j_r j}h_{i_1j_1}\cdots
h_{i_{r-1}j_{r-1}}h_{i_rj_r}
$$
with ${\rm div}P_r=0$. For detail, see \cite{Reilly1973,Barbosa1997,Caoli2007} and the references therein. Then, we have

\begin{thm}\label{1-thm6}
Let $M$ be an $n$-dimensional closed hypersurface in a space form $\mathbb{R}^{n+1}(c)$ and $V$ be a positive twice differentiable function on $M$. If $V{\rm \widehat{Ric}}^{V}_{\phi,\infty}\geq-(n-1)K_1$, where $K_1$ is a nonnegative constant, then for $2\leq r\leq n$,
\begin{align}\label{6-th-111}
\frac{1}{n^2}\int_M V(S_r-\overline{S_r}^V)^2\,d\mu
\leq C_{n,K_1,K_2,\eta_1}\int_M V\Big(|P_r|^2-\frac{(n-r)^2}{n}S_r^2\Big)\,d\mu,
\end{align}
where $\overline{S_r}^V=\frac{\int_M VS_r\,d\mu}{\int_M V\,d\mu}$ and $K_2=\max(V|\nabla \phi|^2)$ for some nonnegative constant $K_2$,
\end{thm}

On the other hand, by applying Theorem \ref{1-thm1}, we also achieve almost Schur lemmas on closed Riemannian manifolds:

\begin{thm}\label{1-thm7}
Let $M$ be an $n$-dimensional closed Riemannian manifold and $V$ be a positive twice differentiable function on $M$.
If $V{\rm \widehat{Ric}}^{V}_{\phi,\infty}\geq-(n-1)K_1$, where $K_1$ is a nonnegative constant, then
\begin{align}\label{7-th-111}
\frac{(n-2)^2}{4n^2}\int_M V(R-\overline{R}^V)^2\,d\mu
\leq C_{n,K_1,K_2,\eta_1}\int_M V\Big(|{\rm Ric}|^2-\frac{1}{n}R^2\Big)\,d\mu,
\end{align}
where
$\overline{R}^V=\frac{\int_M VR\,d\mu}{\int_M V\,d\mu}$ and $K_2=\max(V|\nabla \phi|^2)$ for some nonnegative constant $K_2$.
\end{thm}

The $r$-scalar curvature $\sigma_r$ on a Riemannian manifold $M$ is related to the $k$-th Newton tensor ${T^{(k)}}$ is defined by
$${T^{(k)}}_{i}^{j}=\frac{1}{k!}\sum \delta^{j_1\cdots j_kj}_{i_1\cdots i_ki}A_{i_1}^{j_1}\cdots A_{i_k}^{j_k},$$
where $A$ is the well-known Schouten tensor given by
$$A_{ij}=\frac{1}{n-2}\Big(R_{ij}-\frac{1}{2(n-1)} R\,g_{ij}\Big).$$
Moreover, ${\rm tr_g}(T^{(k)})=(n-k)\sigma_{k}$.
When the metric is locally conformally flat, then we have ${\rm div}{T^{(k)}}=0$, for example, see \cite{Viaclovsky2000,HuLi2008,HuLi2004,GuoHanLi2011}. Therefore, we have

\begin{thm}\label{1-thm8}
Let $M$ be an $n$-dimensional closed locally conformally flat Riemannian manifold and $V$ be a positive twice differentiable function on $M$. If $V{\rm \widehat{Ric}}^{V}_{\phi,\infty}\geq-(n-1)K_1$, where $K_1$ is a nonnegative constant, then for $2\leq k\leq n$,
\begin{align}\label{8-th-111}
\frac{1}{n^2}\int_M V(\sigma_{k}-\overline{\sigma_{k}}^V)^2\,d\mu
\leq C_{n,K_1,K_2,\eta_1}\int_M V\Big(|T^{(k)}|^2-\frac{(n-k)^2}{n}\sigma_{k}^2\Big)\,d\mu,
\end{align}
where $\overline{\sigma_{k}}^V=\frac{\int_M V\sigma_{k}\,d\mu}{\int_M V\,d\mu}$ and $K_2=\max(V|\nabla \phi|^2)$ for some nonnegative constant $K_2$.
\end{thm}

\begin{rem}
From our Theorem \ref{1-thm5} and Theorem \ref{1-thm7}, we can deduce Theorem 1.9 and Theorem 6.1 of Li and Xia in \cite{lixia2019}. Moreover, when $V=1$ and $\phi$ is a constant, our Theorem \ref{1-thm6} and Theorem \ref{1-thm8} become Theorem 1.10 and  Theorem 1.11 of Cheng in \cite{Cheng2014}, respectively.

\end{rem}

In the following, we consider the eigenvalue problem for Wentzell boundary:
\begin{equation}\label{Eig-Int-1}
\Delta_\phi f-\frac{\Delta_\phi V}{V}f=0\  {\rm on }\ M,\ \ \ -\beta\Big(\overline{\Delta}_\phi f-\frac{\overline{\Delta}_\phi V}{V} f\Big)+V\Big(\frac{f}{V}\Big)_\nu=\lambda \frac{f}{V} \ {\rm on }\ \partial M,
\end{equation}
where $\beta$ is a given real constant. When $\beta=0$, the problem \eqref{Eig-Int-1} becomes the following second order Steklov problem:
\begin{equation}\label{Eig-Int-2}
\Delta_\phi f-\frac{\Delta_\phi V}{V}f=0\  {\rm on }\ M,\ \ \ V\Big(\frac{f}{V}\Big)_\nu=p \frac{f}{V} \ {\rm on }\ \partial M.
\end{equation}
Denote by $\lambda_{1,\beta}$ and $p_1$ the first nonzero eigenvalues of \eqref{Eig-Int-1} and \eqref{Eig-Int-2}, respectively. Then, we obtain the following:

\begin{thm}\label{1-thm9}
Let $V$ be a positive twice differentiable function on a given compact Riemannian manifold $M$ with the boundary $\partial M$. Suppose that $VII^V\geq c_1$ and $H_\phi\geq c_2$ for two positive constants $c_1,c_2$. Then we have the following:

(1) If ${\rm \widehat{Ric}}^{V}_{\phi,m}\geq0$, where $m\in (-\infty,0)\cup [n,+\infty)$, then the first nonzero eigenvalue $\eta_1$ of the eigenvalue problem $e^{\phi}{\rm div}(e^{-\phi}V^2\overline{\nabla}\frac{z}{V})=-\eta z$ on the boundary satisfies
\begin{align}\label{9-th-111}
\eta_1\geq c_1c_2,
\end{align}
with the equality holding if and only if $m=n$.

(2) If $V{\rm \widehat{Ric}}^{V}_{\phi,m}\geq-(m-1)K$ for some nonnegative constant $K$, where $m\in (-\infty,0)\cup [n,+\infty)$, then the first nonzero eigenvalue $\lambda_{1,\beta}$ of the eigenvalue problem \eqref{Eig-Int-1} satisfy
\begin{align}\label{9-th-222}
\lambda_{1,\beta}\leq\beta\eta_1+\frac{1}{2c_2}\Big([2\eta_1+(m-1)K]+\sqrt{[2\eta_1+(m-1)K]^2-4c_1c_2\eta_1}\Big),
\end{align}
with the equality holding if and only if $m=n$.

\end{thm}

\begin{thm}\label{1-thm10}
Let $V$ be a positive twice differentiable function on a given compact Riemannian manifold $M$ with the boundary $\partial M$. Suppose that $VII^V\geq  c_1$ and $H_\phi\geq c_2$ for two positive constants $c_1,c_2$. Then we have the following:

(1) If $V{\rm \widehat{Ric}}^{V}_{\phi,m}\geq-(m-1)K$ for some nonnegative constant $K$, where $m\in (-\infty,0)\cup [n,+\infty)$, then the first nonzero eigenvalue $p_1$ of the eigenvalue problem \eqref{Eig-Int-2} satisfy
\begin{align}\label{10-th-111}
p_1>\frac{c_1\eta_1}{2\eta_1+(m-1)K}.
\end{align}

(2) If ${\rm \widehat{Ric}}^{V}_{\phi,m}\geq0$, where $m\in (-\infty,0)\cup [n,+\infty)$, then the first nonzero eigenvalue $\lambda_{1,\beta}$ of the eigenvalue problem \eqref{Eig-Int-1} satisfy
\begin{align}\label{10-th-222}
\lambda_{1,\beta}\geq\frac{c_1}{2}\Big(1+\beta c_2+\sqrt{\beta^2 c_2^2+2\beta c_2}\Big).
\end{align}

\end{thm}

\begin{rem}
When $V=1$, the formula \eqref{9-th-111} is corresponding to Theorem 1.6 of Huang and Ruan \cite{HuangRuan2014} and
the estimate \eqref{9-th-222} is corresponding to Theorem 1.1 of Wang and Xia in \cite{Wang2018}. Our Theorem \ref{1-thm10} generalizes Theorem 1.3 of Wang and Xia in \cite{Wang2018}. Moreover, by letting $K=0$, $V=1$ and $\beta\rightarrow 0$ in \eqref{9-th-222}, we obtain the estimate (1.4) in Theorem 1.1 of Wang and Xia in \cite{Wang2009}.
\end{rem}

\begin{ack}
We would like to thank Dr. Fanqi Zeng for helpful discussions which make the paper more readable.
\end{ack}

\section{Proof of results}
\subsection{Proof of Theorem \ref{1-thm1}}
Using the divergence theorem, we have
\begin{align}\label{2-Sec-Prof-Th-1}
\int_MVf_{ij}^2\,d\mu=&\int_{\partial M}\frac{1}{2}V|\nabla f|^2_{\nu}\,d\sigma-\int_Me^\phi(e^{-\phi}Vf_{ij})_jf_i\,d\mu\notag\\
=&\int_{\partial M}\frac{1}{2}V|\nabla f|^2_{\nu}\,d\sigma-\int_M(f_{ij}V_j+V f_{ij,j}-Vf_{ij}\phi_j)f_i\,d\mu\notag\\
=&\int_{\partial M}\frac{1}{2}V|\nabla f|^2_{\nu}\,d\sigma-\int_M\Big[\frac{1}{2}|\nabla f|^2_iV_i+V(\Delta_\phi f)_i f_i
+V{\rm Ric}_{\phi}(\nabla f,\nabla f)\Big]\,d\mu\notag\\
=&\int_{\partial M}\Big[\frac{1}{2}V|\nabla f|^2_{\nu}-\frac{1}{2}|\nabla f|^2V_{\nu}-V(\Delta_\phi f)f_{\nu}\Big]\,d\sigma+\int_M\Big[\frac{1}{2}|\nabla f|^2\Delta_\phi V\notag\\
&+(\Delta_\phi f)V_if_i+V(\Delta_\phi f)^2-V{\rm Ric}_{\phi}(\nabla f,\nabla f)\Big]\,d\mu,
\end{align}
\begin{align}\label{2-Sec-Prof-Th-2}
-2\int_Mff_{ij}V_{ij}\,d\mu=&-2\int_{\partial M}ff_{i\nu}V_{i}\,d\sigma+2\int_M[f_{ij}V_if_j+f(\Delta_\phi f)_i V_i\notag\\
&+f{\rm Ric}_{\phi}(\nabla f,\nabla V)]\,d\mu\notag\\
=&\int_{\partial M}[2f(\Delta_\phi f)V_\nu-2ff_{i\nu}V_{i}+|\nabla f|^2V_\nu]\,d\sigma\notag\\
&+\int_M[-|\nabla f|^2\Delta_\phi V-2f(\Delta_\phi f)(\Delta_\phi V)\notag\\
&-2(\Delta_\phi f) f_i V_i+2f{\rm Ric}_{\phi}(\nabla f,\nabla V)]\,d\mu
\end{align}
and
\begin{align}\label{2-Sec-Prof-Th-3}
\int_M\frac{f^2}{V}V_{ij}^2\,d\mu=&\frac{1}{2}\int_{\partial M}\frac{f^2}{V}|\nabla V|^2_{\nu}\,d\sigma-\int_M\Big[\Big(\frac{f^2}{V}\Big)_jV_{ij}V_i\notag\\
&+\frac{f^2}{V}(\Delta_\phi V)_iV_i+\frac{f^2}{V}{\rm Ric}_{\phi}(\nabla V,\nabla V)\Big]\,d\mu\notag\\
=&\int_{\partial M}\Big[\frac{1}{2}\frac{f^2}{V}|\nabla V|^2_{\nu}-\frac{f^2}{V}(\Delta_\phi V)V_\nu\Big]\,d\sigma\notag\\
&+\int_M\Big[\frac{f^2}{V}(\Delta_\phi V)^2+(\Delta_\phi V)\Big(\frac{f^2}{V}\Big)_iV_i-\Big(\frac{f^2}{V}\Big)_jV_{ij}V_i\notag\\
&-\frac{f^2}{V}{\rm Ric}_{\phi}(\nabla V,\nabla V)\Big]\,d\mu.
\end{align}
Let $A_{ij}=f_{ij}-\frac{f}{V} V_{ij}$. Then we have $A_{ij}^2=f_{ij}^2-2\frac{f}{V}f_{ij}V_{ij}+\frac{f^2}{V^2}V_{ij}^2$ and ${\rm tr}(A_{ij})=\Delta f-\frac{f}{V}\Delta V$. It follows from \eqref{2-Sec-Prof-Th-1}-\eqref{2-Sec-Prof-Th-3} that
\begin{align}\label{2-Sec-Prof-Th-4}
\int_MVA_{ij}^2\,d\mu=&\int_{\partial M}\Big[\frac{1}{2}V|\nabla f|^2_{\nu}+\frac{1}{2}|\nabla f|^2V_{\nu}-V(\Delta_\phi f)f_{\nu}+2f(\Delta_\phi f)V_\nu\notag\\
&-2ff_{i\nu}V_{i}+\frac{1}{2}\frac{f^2}{V}|\nabla V|^2_{\nu}-\frac{f^2}{V}(\Delta_\phi V)V_\nu\Big]\,d\sigma\notag\\
&+\int_M\Big[-\frac{1}{2}|\nabla f|^2\Delta_\phi V-(\Delta_\phi f)V_if_i+(\Delta_\phi V)\Big(\frac{f^2}{V}\Big)_iV_i\notag\\
&-\Big(\frac{f^2}{V}\Big)_jV_{ij}V_i\Big]\,d\mu
+\int_M V\Big[\Big(\Delta_\phi f-\frac{f}{V}\Delta_\phi V\Big)^2\notag\\
&-{\rm Ric}_{\phi}\Big(V\nabla\frac{f}{V},V\nabla\frac{f}{V}\Big)\Big]\,d\mu,
\end{align}
where we note that $\nabla f-\frac{f}{V}\nabla V=V\nabla\frac{f}{V}$. It is easy to check that
$$\aligned
{}[-(\Delta_\phi V)g_{ij}+V_{ij}]\Big(V\nabla\frac{f}{V},V\nabla\frac{f}{V}\Big)=&-(\Delta_\phi V)\Big[|\nabla f|^2-\Big(\frac{f^2}{V}\Big)_iV_{i}\Big]\\
&+V_{ij}f_if_j-\Big(\frac{f^2}{V}\Big)_jV_{ij}V_i,
\endaligned$$
which is equivalent to
\begin{align}\label{2-Sec-Prof-Th-7}
(\Delta_\phi V)\Big(\frac{f^2}{V}\Big)_iV_i-\Big(\frac{f^2}{V}\Big)_jV_{ij}V_i=&[-(\Delta_\phi V)g_{ij}+V_{ij}]\Big(V\nabla\frac{f}{V},V\nabla\frac{f}{V}\Big)\notag\\
&+(\Delta_\phi V)|\nabla f|^2-V_{ij}f_if_j.
\end{align}
Therefore, \eqref{2-Sec-Prof-Th-4} becomes
\begin{align}\label{2-Sec-Prof-Th-8}
\int_MVA_{ij}^2\,d\mu=&\int_{\partial M}\Big[\frac{1}{2}V|\nabla f|^2_{\nu}+\frac{1}{2}|\nabla f|^2V_{\nu}-V(\Delta_\phi f)f_{\nu}+2f(\Delta_\phi f)V_\nu\notag\\
&-2ff_{i\nu}V_{i}+\frac{1}{2}\frac{f^2}{V}|\nabla V|^2_{\nu}-\frac{f^2}{V}(\Delta_\phi V)V_\nu\Big]\,d\sigma\notag\\
&+\int_M\Big[\frac{1}{2}|\nabla f|^2\Delta_\phi V-(\Delta_\phi f)V_if_i-V_{ij}f_if_j\Big]\,d\mu\notag\\
&+\int_M V\Big[\Big(\Delta_\phi f-\frac{f}{V}\Delta_\phi V\Big)^2-{\rm \widehat{Ric}}^{V}_{\phi,\infty}\Big(V\nabla\frac{f}{V},V\nabla\frac{f}{V}\Big)\Big]\,d\mu\notag\\
=&\int_{\partial M}\Big[\frac{1}{2}V|\nabla f|^2_{\nu}+|\nabla f|^2V_{\nu}-V(\Delta_\phi f)f_{\nu}+2f(\Delta_\phi f)V_\nu\notag\\
&-2ff_{i\nu}V_{i}+\frac{1}{2}\frac{f^2}{V}|\nabla V|^2_{\nu}-\frac{f^2}{V}(\Delta_\phi V)V_\nu-f_\nu\nabla f\nabla V \Big]\,d\sigma\notag\\
&+\int_M V\Big[\Big(\Delta_\phi f-\frac{f}{V}\Delta_\phi V\Big)^2-{\rm \widehat{Ric}}^{V}_{\phi,\infty}\Big(V\nabla\frac{f}{V},V\nabla\frac{f}{V}\Big)\Big]\,d\mu
\end{align}
from
$$\aligned
-\int_MV_{ij}f_if_j\,d\mu=&-\int_{\partial M}\nabla f\nabla V f_\nu\,d\sigma+\int_M[V_if_{ij}f_j+(\Delta_\phi f)V_if_i]\,d\mu\\
=&\int_{\partial M}\Big(\frac{1}{2}|\nabla f|^2V_{\nu}-f_\nu\nabla f\nabla V\Big)\,d\sigma\\
&+\int_M\Big(-\frac{1}{2}|\nabla f|^2\Delta_\phi V+(\Delta_\phi f)V_if_i\Big)\,d\mu.
\endaligned$$

Let $\{e_i\}_{i=1}^n$ be an orthonormal frame field along $\partial M$ such that $\{e_\alpha\}_{\alpha=1}^{n-1}$ are
tangent to $\partial M$ and $e_{n}=\nu$ is normal to $\partial M$. Then we have the following:
$$|\nabla f|^2_\nu=2\sum_{\alpha=1}^{n-1}f_\alpha f_{\alpha\nu}+2f_{\nu\nu}f_{\nu}=2\overline{\nabla} f\overline{\nabla} f_\nu-2II(\overline{\nabla}f,\overline{\nabla} f)+2f_{\nu\nu}f_{\nu},$$
$$|\nabla f|^2=|\overline{\nabla} f|^2+f_\nu^2,$$
$$\Delta_\phi f=\overline{\Delta}_\phi f+H_\phi f_{\nu}+f_{\nu\nu},$$
$$f_{i\nu}V_{i}=f_{\alpha\nu}V_{\alpha}+f_{\nu\nu}V_{\nu}=\overline{\nabla} V\overline{\nabla} f_\nu-II(\overline{\nabla}f,\overline{\nabla} V)+f_{\nu\nu}V_{\nu},$$
$$|\nabla V|^2_\nu=2\sum_{\alpha=1}^{n-1}V_\alpha V_{\alpha\nu}+2V_{\nu\nu}V_{\nu}=2\overline{\nabla} V\overline{\nabla} V_\nu-2II(\overline{\nabla}V,\overline{\nabla} V)+2V_{\nu\nu}V_{\nu},$$
$$\Delta_\phi V=\overline{\Delta}_\phi V+H_\phi V_{\nu}+V_{\nu\nu}$$
and
$$f_\nu\nabla f\nabla V=f_\nu\overline{\nabla} f\overline{\nabla}V+V_{\nu}f_{\nu}^2.$$
Therefore, on $\partial M$,
\begin{align}\label{2-Sec-Prof-Th-9}
\frac{1}{2}&V|\nabla f|^2_{\nu}+|\nabla f|^2V_{\nu}-V(\Delta_\phi f)f_{\nu}+2f(\Delta_\phi f)V_\nu\notag\\
&-2ff_{i\nu}V_{i}+\frac{1}{2}\frac{f^2}{V}|\nabla V|^2_{\nu}-\frac{f^2}{V}(\Delta_\phi V)V_\nu- f_\nu\nabla f\nabla V\notag\\
=&-V(II-(\ln V)_\nu g)\Big(\overline{\nabla}f-\frac{f}{V}\overline{\nabla}V,\overline{\nabla}f-\frac{f}{V}\overline{\nabla}V\Big)\notag\\
&-VH_\phi\Big(f_{\nu}-\frac{f}{V}V_{\nu}\Big)^2+2\frac{f}{V}V_{\nu}\overline{\nabla} f\overline{\nabla}V-\frac{f^2}{V^2}V_{\nu}|\overline{\nabla}V|^2\notag\\
&+V\overline{\nabla}f\overline{\nabla}f_{\nu}-Vf_{\nu}\overline{\Delta}_\phi f+2fV_{\nu}\overline{\Delta}_\phi f
-2f\overline{\nabla}V\overline{\nabla}f_{\nu}\notag\\
&+\frac{f^2}{V}\overline{\nabla}V\overline{\nabla}V_{\nu}-\frac{f^2}{V}V_{\nu}\overline{\Delta}_\phi V-f_\nu\overline{\nabla} f\overline{\nabla}V.
\end{align}
Using the fact that the boundary $\partial M$ is closed, then we have
\begin{align}\label{2-Sec-Prof-Th-10}
&\int_{\partial M}\Big(V\overline{\nabla}f\overline{\nabla}f_{\nu}-2f\overline{\nabla}V\overline{\nabla}f_{\nu}
+\frac{f^2}{V}\overline{\nabla}V\overline{\nabla}V_{\nu}\Big)\,d\sigma\notag\\
=&\int_{\partial M}\Big(-f_{\nu}\overline{\nabla}V\overline{\nabla}f-Vf_{\nu}\overline{\Delta}_\phi f+2 f_{\nu}\overline{\nabla}f\overline{\nabla}V+2 f f_{\nu}\overline{\Delta}_\phi V\notag\\
&-V_{\nu}\overline{\nabla}\frac{f^2}{V}\overline{\nabla}V-V_{\nu}\frac{f^2}{V}\overline{\Delta}_\phi V\Big)\,d\sigma
\end{align}
and
\begin{align}\label{2-Sec-Prof-Th-11}
&\int_{\partial M}\Big[\frac{1}{2}V|\nabla f|^2_{\nu}+|\nabla f|^2V_{\nu}-V(\Delta_\phi f)f_{\nu}+2f(\Delta_\phi f)V_\nu\notag\\
&-2ff_{i\nu}V_{i}+\frac{1}{2}\frac{f^2}{V}|\nabla V|^2_{\nu}-\frac{f^2}{V}(\Delta_\phi V)V_\nu- f_\nu\nabla f\nabla V\Big]\,d\sigma\notag\\
=&\int_{\partial M}\Big[-VII^V\Big(V\overline{\nabla}\frac{f}{V},V\overline{\nabla}\frac{f}{V}\Big)-VH_\phi\Big(f_{\nu}-(\ln V)_\nu f \Big)^2\notag\\
&-2V^2\Big(\frac{f}{V}\Big)_\nu\Big(\overline{\Delta}_\phi f-\frac{\overline{\Delta}_\phi V}{V} f\Big)\Big]\,d\sigma,
\end{align}
which concludes the proof of Theorem \ref{1-thm1}.

\subsection{Proof of Theorem \ref{1-thm3}}
We consider the following Dirichlet boundary problem:
$$\Delta_\phi f-\frac{\Delta_\phi V}{V}f=1 \  {\rm on }\ M,\ \ \ f=0\ {\rm on }\ \partial M.$$
Then from \eqref{2-corr-111}, we obtain
\begin{align}\label{2-Sec-Prof-Th2-1}
\int_{\partial M}VH_\phi f_{\nu}^2\,d\sigma\leq\frac{m-1}{m}\int_MV\,d\mu.
\end{align}
On the other hand,
\begin{align}\label{2-Sec-Prof-Th2-2}
\int_MV\,d\mu=&\int_M(V\Delta_\phi f-f \Delta_\phi V)\,d\mu\notag\\
=&\int_{\partial M}Vf_{\nu}\,d\sigma\notag\\
\leq&\Big(\int_{\partial M}\frac{V}{H_\phi}\,d\sigma\Big)^{\frac{1}{2}}\Big(\int_{\partial M}VH_\phi f_{\nu}^2\,d\sigma\Big)^{\frac{1}{2}},
\end{align}
which together with \eqref{2-Sec-Prof-Th2-1} gives the formula \eqref{3-th-111}.

If the equality in \eqref{3-th-111} is attained, then from $\Delta_\phi f-\frac{\Delta_\phi V}{V}f=1$ and \eqref{Int-5}, we have $\Delta f-\frac{\Delta  V}{V}f=\frac{n}{m}$ and $V\nabla \phi\nabla \frac{f}{V}=-\frac{m-n}{m}$. If $m>n$, then \eqref{Int-5} shows that
$$V\Delta f-f\Delta V+\frac{n}{m-n}V^2\nabla \phi\nabla \frac{f}{V}=0,$$
which is equivalent to
\begin{align}\label{2-Sec-Prof-Th2-3}
e^{-\frac{n}{m-n}\phi}\Big[e^{\frac{n}{m-n}\phi}V^2\Big(\frac{f}{V}\Big)_i\Big]_i=&\Big[V^2\Big(\frac{f}{V}\Big)_i\Big]_i
+\frac{n}{m-n}V^2\nabla \phi\nabla \frac{f}{V}\notag\\
=&0,
\end{align}
where we notice that $V\Delta \frac{f}{V}=\Delta f-\frac{f}{V}\Delta V-2\nabla V\nabla \frac{f}{V}$.
Multiplying both sides of \eqref{2-Sec-Prof-Th2-3} with $\frac{f}{V}e^{\frac{n}{m-n}\phi}$ gives
\begin{align}\label{2-Sec-Prof-Th2-4}
0=&\int_M\frac{f}{V}\Big[e^{\frac{n}{m-n}\phi}V^2\Big(\frac{f}{V}\Big)_i\Big]_i\,dv_g\notag\\
=&-\int_M V^2\Big|\nabla\frac{f}{V}\Big|^2 e^{\frac{n}{m-n}\phi}\,dv_g,
\end{align}
and then $f=\theta V$, where $\theta$ is a constant. This contradicts with $\Delta_\phi f-\frac{\Delta_\phi V}{V}f=1$. Therefore, we have that if the equality in \eqref{3-th-111} is attained, then $m=n$ and $\phi$ must be constant. Using the similar conclusions as in Theorem 1.3 of \cite{lixia2019}, we complete the proof.

\subsection{Proof of Theorem \ref{1-thm4}}
We consider the following Dirichlet boundary problem:
$$\Delta_\phi f-\frac{\Delta_\phi V}{V}f=1 \  {\rm on }\ M,\ \ \ V\Big(\frac{f}{V}\Big)_\nu=c\ {\rm on }\ \partial M,$$
where
$$c=\frac{\int_MV\,d\mu}{\int_{\partial M}V\,d\sigma}.$$
The existence and uniqueness of the solution to above equation is due to Fredholm alternative.
Then from \eqref{2-corr-111}, we obtain
\begin{align}\label{4-th-Proof-1}
0\leq&\int_{\partial M}[-c^2VH_\phi-2c(V\overline{\Delta}_\phi f-f\overline{\Delta}_\phi V)]\,d\sigma+\frac{m-1}{m}\int_MV\,d\mu\notag\\
=&-c^2\int_{\partial M}VH_\phi\,d\sigma+\frac{m-1}{m}\int_MV\,d\mu,
\end{align}
which gives \eqref{4-th-222}.

Similarly, if the equality in \eqref{4-th-222} is attained the inequality in \eqref{4-th-111} is strict, then we have that $m=n$ and $\phi$ is a constant. Then, according to the arguments as in Theorem 1.7 of \cite{lixia2019} finishes the proof.

\subsection{Proof of Theorems \ref{1-thm5}-\ref{1-thm8}}
Firstly, we prove the following

\begin{prop}\label{5-prop5}
Let $V$ be a positive twice differentiable function on an $n$-dimensional closed Riemannian manifold $M$ with $V{\rm \widehat{Ric}}^{V}_{\phi,\infty}\geq-(n-1)K_1$, where $K_1$ is a nonnegative constant. If the symmetric $(2,0)$-tensor field $T$ defined on $M$ satisfies ${\rm div}T=c\nabla({\rm tr}T)$, where $c$ is a constant, then
\begin{align}\label{5-prop-111}
\frac{(nc-1)^2}{n^2}\int_M V({\rm tr}T-\overline{{\rm tr}T}^V)^2\,d\mu
\leq C_{n,K_1,K_2,\eta_1}\int_M V\Big|T-\frac{1}{n}({\rm tr}T)g\Big|^2\,d\mu,
\end{align}
where $K_2=\max(V|\nabla \phi|^2)$ for some nonnegative constant $K_2$,
$$\overline{{\rm tr}T}^V=\frac{\int_M V{\rm tr}T\,d\mu}{\int_M V\,d\mu}$$
and the constant $C_{n,K_1,K_2,\eta_1}$ is given by \eqref{const-1}. Here $\eta_1$ is the first nonzero eigenvalue of problem $\Delta_\phi f-\frac{\Delta_\phi V}{V}f=-\eta \frac{f}{V}$.
\end{prop}

\proof
Let $f$ be a solution to the following problem:
$$\Delta_\phi f-\frac{\Delta_\phi V}{V}f={\rm tr}T-\overline{{\rm tr}T}^V \  {\rm on }\ M.$$
Denote by $A_{ij}=f_{ij}-\frac{f}{V}V_{ij}$ and $\mathring{A}_{ij}=A_{ij}-\frac{{\rm tr}(A_{ij})}{n}g_{ij}$ with ${\rm tr}(A_{ij})=\Delta f-\frac{f}{V}\Delta V$. Then, we have ${\rm tr}(\mathring{A}_{ij})=0$ and
\begin{align}\label{5-th-Proof-1}
&\int_M V({\rm tr}T-\overline{{\rm tr}T}^V)^2\,d\mu\notag\\
=&\int_M ({\rm tr}T-\overline{{\rm tr}T}^V)(V\Delta_\phi f-f\Delta_\phi V)\,d\mu\notag\\
=&-\int_M (V\nabla f-f\nabla V)\nabla {\rm tr}T\,d\mu.
\end{align}
Using $\nabla{\rm tr}T=\frac{1}{c}{\rm div}T$, we have $\frac{nc-1}{n}\nabla{\rm tr}T={\rm div}\mathring{T}$, where $\mathring{T}_{ij}=T_{ij}-\frac{1}{n} ({\rm tr}T)g_{ij}$. Then \eqref{5-th-Proof-1} gives
\begin{align}\label{5-th-Proof-2}
&\int_M V({\rm tr}T-\overline{{\rm tr}T}^V)^2\,d\mu\notag\\
=&-\frac{n}{nc-1}\int_M \mathring{T}_{ij,j}(Vf_i-fV_i)\,d\mu\notag\\
=&\frac{n}{nc-1}\int_M V\mathring{T}_{ij}B_{ij}\,d\mu\notag\\
=&\frac{n}{nc-1}\int_M V\mathring{T}_{ij}\mathring{B}_{ij}\,d\mu\notag\\
\leq&\frac{n}{|nc-1|}\Big(\int_M V|\mathring{T}_{ij}|^2\,d\mu\Big)^{\frac{1}{2}}\Big(\int_M V|\mathring{B}_{ij}|^2\,d\mu\Big)^{\frac{1}{2}},
\end{align}
where $\mathring{B}_{ij}=B_{ij}-\frac{1}{n} ({\rm tr}B)g_{ij}$ and
$$B_{ij}=A_{ij}-\frac{1}{2}V\Big[\Big(\frac{f}{V}\Big)_i\phi_j+\phi_i\Big(\frac{f}{V}\Big)_j\Big].$$
It is easy to check that ${\rm tr}B=\Delta_\phi f-\frac{\Delta_\phi V}{V}f$.
By virtue of the Cauchy inequality, we have
\begin{align}\label{5-th-Proof-3}
|\mathring{B}_{ij}|^2=&|B_{ij}|^2-\frac{1}{n}\Big(\Delta_\phi f-\frac{\Delta_\phi V}{V}f\Big)^2\notag\\
\leq&(1+\delta)|A_{ij}|^2+\Big(1+\frac{1}{\delta}\Big)V^2\Big|\nabla\frac{f}{V}\Big|^2|\nabla\phi|^2\notag\\
&-\frac{1}{n}\Big(\Delta_\phi f-\frac{\Delta_\phi V}{V}f\Big)^2,
\end{align}
where $\delta$ is a positive constant to be determined, and then
\begin{align}\label{5-th-Proof-4}
\int_M V|\mathring{B}_{ij}|^2\,d\mu\leq&(1+\delta)\int_MV|A_{ij}|^2\,d\mu+\Big(1+\frac{1}{\delta}\Big)\int_MV^3\Big|\nabla\frac{f}{V}\Big|^2|\nabla\phi|^2\,d\mu\notag\\
&-\frac{1}{n}\int_MV\Big(\Delta_\phi f-\frac{\Delta_\phi V}{V}f\Big)^2\,d\mu.
\end{align}
On other hand, applying \eqref{1-th-111} on closed manifold $M$, we obtain
\begin{align}\label{5-th-Proof-5}
\int_MV|A_{ij}|^2\,d\mu=&\int_MV\Big|\nabla^2f-\frac{f}{V} \nabla^2V\Big|^2\,d\mu\notag\\
=&\int_M\Big[V\Big(\Delta_\phi f-\frac{f}{V}\Delta_\phi V\Big)^2-V{\rm \widehat{Ric}}^{V}_{\phi,\infty}\Big(V\nabla\frac{f}{V},V\nabla\frac{f}{V}\Big)\Big]\,d\mu.
\end{align}
Thus, \eqref{5-th-Proof-4} becomes
\begin{align}\label{5-th-Proof-6}
\int_M V|\mathring{B}_{ij}|^2\,d\mu\leq&\Big(\frac{n-1}{n}+\delta\Big)\int_MV\Big(\Delta_\phi f-\frac{\Delta_\phi V}{V}f\Big)^2\,d\mu+\Big[\Big(1+\frac{1}{\delta}\Big)K_2\notag\\
&+(n-1)(1+\delta)K_1\Big]\int_MV^2\Big|\nabla\frac{f}{V}\Big|^2\,d\mu.
\end{align}
Using the Rayleigh-Reitz principle, we have the first nonzero eigenvalue $\eta_1$ of $\Delta_\phi f-\frac{\Delta_\phi V}{V}f=-\eta \frac{f}{V}$ can also be characterized by
\begin{align}\label{5-th-Proof-7}
\eta_1=\inf\limits_{f\in C^\infty(M)} \frac{\int_MV^2|\nabla\frac{f}{V}|^2\,d\mu}{\int_MV(\frac{f}{V})^2\,d\mu}.
\end{align}
Then,
\begin{align}\label{5-th-Proof-8}
\int_MV^2\Big|\nabla\frac{f}{V}\Big|^2\,d\mu=&-\int_MV\Big(\frac{f}{V}\Big)\Big(\Delta_\phi f-\frac{\Delta_\phi V}{V}f\Big)\,d\mu\notag\\
\leq&\Big[\int_MV\Big(\frac{f}{V}\Big)^2\,d\mu\Big]^{\frac{1}{2}}\Big[\int_MV\Big(\Delta_\phi f-\frac{\Delta_\phi V}{V}f\Big)^2\,d\mu\Big]^{\frac{1}{2}}\notag\\
\leq&\Big[\frac{1}{\eta_1}\int_MV^2\Big|\nabla\frac{f}{V}\Big|^2\,d\mu\Big]^{\frac{1}{2}}\Big[\int_MV\Big(\Delta_\phi f-\frac{\Delta_\phi V}{V}f\Big)^2\,d\mu\Big]^{\frac{1}{2}},
\end{align}
which gives
\begin{align}\label{5-th-Proof-9}
\int_MV^2\Big|\nabla\frac{f}{V}\Big|^2\,d\mu\leq&\frac{1}{\eta_1}\int_MV\Big(\Delta_\phi f-\frac{\Delta_\phi V}{V}f\Big)^2\,d\mu\notag\\
=&\frac{1}{\eta_1}\int_MV({\rm tr}T-\overline{{\rm tr}T}^V)^2\,d\mu.
\end{align}
Inserting \eqref{5-th-Proof-9} into \eqref{5-th-Proof-6}, we have
\begin{align}\label{5-th-Proof-10}
\int_M V|\mathring{B}_{ij}|^2\,d\mu\leq&\frac{1}{\eta_1}\Big[\frac{n-1}{n}\eta_1+\delta\eta_1+\Big(1+\frac{1}{\delta}\Big)K_2\notag\\
&+(n-1)(1+\delta)K_1\Big]\int_MV({\rm tr}T-\overline{{\rm tr}T}^V)^2\,d\mu\notag\\
=&\frac{1}{\eta_1}\Big[\frac{n-1}{n}\eta_1+(n-1)K_1+K_2+[\eta_1+(n-1)K_1]\delta\notag\\
&+\frac{K_2}{\delta}\Big]\int_MV({\rm tr}T-\overline{{\rm tr}T}^V)^2\,d\mu.
\end{align}
Minimizing the $\delta$ in \eqref{5-th-Proof-10} by taking
$$\delta=\sqrt{\frac{K_2}{\eta_1+(n-1)K_1}},$$
we obtain
\begin{align}\label{5-th-Proof-11}
\int_M V|\mathring{B}_{ij}|^2\,d\mu
\leq C_{n,K_1,K_2,\eta_1}\int_MV({\rm tr}T-\overline{{\rm tr}T}^V)^2\,d\mu.
\end{align}
Therefore, combining \eqref{5-th-Proof-11} with \eqref{5-th-Proof-2} concludes the proof of \eqref{5-prop-111} and the proof of Proposition \ref{5-prop5} is finished.
\endproof

Since the ambient space is a space form, we have the well-known Codazzi equation:
$$II_{ij,j}=H_i.$$
Hence, we complete the proof of Theorems \ref{1-thm5} and \ref{1-thm6} by taking $c=1$ and $c=0$ in the formula \eqref{5-prop-111}, respectively.

Using the second Bianchi identity, the Ricci curvature $R_{ij}$ is related to the scalar curvature $R$ by $R_{ij,j}=\frac{1}{2}R_{,i}$. Thus, the proof of theorems \ref{1-thm7} and \ref{1-thm8} follows by taking $c=\frac{1}{2}$ and $c=0$ in the formula \eqref{5-prop-111}, respectively.

\subsection{Proof of Theorem \ref{1-thm9}}
(1) Let $f$ be the solution of the following equation:
$$\Delta_\phi f-\frac{\Delta_\phi V}{V}f=0 \  {\rm on }\ M,\ \ \ f=z\ {\rm on }\ \partial M,$$
where $z$ satisfies $\overline{\Delta}_\phi z-\frac{\overline{\Delta}_\phi V}{V} z=-\eta_1 \frac{z}{V}$, that is, $z$ satisfies $e^{\phi}{\rm div}(e^{-\phi}V^2\nabla\frac{z}{V})=-\eta_1z$.
Then from \eqref{2-corr-111}, we achieve
\begin{align}\label{9-thm-Proof-1}
0\geq\int_{\partial M}\Big[c_1V^2\Big|\overline{\nabla}\frac{z}{V}\Big|^2+c_2V^3\Big(\Big(\frac{f}{V}\Big)_\nu\Big)^2
-2\eta_1V^2\Big(\frac{f}{V}\Big)_\nu\frac{z}{V}\Big]\,d\sigma.
\end{align}
Applying
\begin{align}\label{9-thm-Proof-2}
\int_{\partial M}V^2\Big|\overline{\nabla}\frac{z}{V}\Big|^2\,d\sigma=&-\int_{\partial M}\frac{z}{V}(V\overline{\Delta}_\phi z-z\overline{\Delta}_\phi V)\,d\sigma\notag\\
=&\eta_1\int_{\partial M}V\Big(\frac{z}{V}\Big)^2\,d\sigma
\end{align}
in \eqref{9-thm-Proof-1} yields
\begin{align}\label{9-thm-Proof-3}
0\geq&\int_{\partial M}\Big[c_1\eta_1V\Big(\frac{z}{V}\Big)^2+c_2V^3\Big(\Big(\frac{f}{V}\Big)_\nu\Big)^2
-2\eta_1V^2\Big(\frac{f}{V}\Big)_\nu\frac{z}{V}\Big]\,d\sigma\notag\\
\geq&\Big(c_1\eta_1-\frac{\eta_1^2}{c_2}\Big)\int_{\partial M}V\Big(\frac{z}{V}\Big)^2\,d\sigma,
\end{align}
which shows that $\eta_1\geq c_1c_2$ and the estimate \eqref{9-th-111} follows.

(2) Let $f$ be the solution of the following equation:
$$\Delta_\phi f-\frac{\Delta_\phi V}{V}f=0\  {\rm on }\ M,\ \ \ f=z\ {\rm on }\ \partial M,$$
where $z$ satisfies $\overline{\Delta}_\phi z-\frac{\overline{\Delta}_\phi V}{V} z=-\eta_1 \frac{z}{V}$, that is, $z$ satisfies $e^{\phi}{\rm div}(e^{-\phi}V^2\nabla\frac{z}{V})=-\eta_1z$.
Then from \eqref{2-corr-111}, we also have
\begin{align}\label{9-thm-Proof-4}
(m-1)K&\int_{\partial M}V^2\Big(\frac{f}{V}\Big)_\nu\frac{z}{V}\,d\sigma\notag\\
=&(m-1)K\int_MV^2\Big|\nabla\frac{f}{V}\Big|^2\,d\mu\notag\\
\geq&\int_{\partial M}\Big[c_1V^2\Big|\overline{\nabla}\frac{z}{V}\Big|^2+c_2V^3\Big(\Big(\frac{f}{V}\Big)_\nu\Big)^2
-2\eta_1V^2\Big(\frac{f}{V}\Big)_\nu\frac{z}{V}\Big]\,d\sigma\notag\\
=&\int_{\partial M}\Big[c_1\eta_1V\Big(\frac{z}{V}\Big)^2+c_2V^3\Big(\Big(\frac{f}{V}\Big)_\nu\Big)^2
-2\eta_1V^2\Big(\frac{f}{V}\Big)_\nu\frac{z}{V}\Big]\,d\sigma,
\end{align}
which gives
\begin{align}\label{9-thm-Proof-5}
0\geq&\int_{\partial M}\Big[c_1\eta_1V\Big(\frac{z}{V}\Big)^2+c_2V^3\Big(\Big(\frac{f}{V}\Big)_\nu\Big)^2
-[2\eta_1+(m-1)K]V^2\Big(\frac{f}{V}\Big)_\nu\frac{z}{V}\Big]\,d\sigma\notag\\
\geq&c_2\int_{\partial M}V^3\Big(\Big(\frac{f}{V}\Big)_\nu\Big)^2\,d\sigma+c_1\eta_1\int_{\partial M}V\Big(\frac{z}{V}\Big)^2\,d\sigma\notag\\
&-[2\eta_1+(m-1)K]\Big(\int_{\partial M}V\Big(\frac{z}{V}\Big)^2\,d\sigma\Big)^{\frac{1}{2}}\Big(\int_{\partial M}V^3\Big(\Big(\frac{f}{V}\Big)_\nu\Big)^2\,d\sigma\Big)^{\frac{1}{2}}.
\end{align}
Therefore, we have proved that
\begin{align}\label{9-thm-Proof-6}
c_2x^2-[2\eta_1+(m-1)K]x+c_1\eta_1\leq0,
\end{align}
where
$$x=\Big[\int_{\partial M}V^3\Big(\Big(\frac{f}{V}\Big)_\nu\Big)^2\,d\sigma\Big/\int_{\partial M}V\Big(\frac{z}{V}\Big)^2\,d\sigma\Big]^{\frac{1}{2}}.$$
Solving this quadratic inequality with respect to $x$ gives
\begin{align}\label{9-thm-Proof-7}
x\leq\frac{1}{2c_2}\Big([2\eta_1+(m-1)K]+\sqrt{[2\eta_1+(m-1)K]^2-4c_1c_2\eta_1}\Big).
\end{align}
On the other hand, from the Rayleigh-Ritz formula(we may refer to \cite{Kuttler1968,Wang2018}), we have
\begin{align}\label{9-thm-Proof-8}
\lambda_{1,\beta}\leq&\beta\eta_1+\frac{\int_MV^2\Big|\nabla\frac{f}{V}\Big|^2\,d\mu}{\int_{\partial M}V(\frac{z}{V})^2\,d\sigma}\notag\\
=&\beta\eta_1+\frac{\int_{\partial M}V^2(\frac{f}{V})_\nu\frac{z}{V}\,d\sigma}{\int_{\partial M}V(\frac{z}{V})^2\,d\sigma}\notag\\
\leq&\beta\eta_1+\Big[\int_{\partial M}V^3\Big(\Big(\frac{f}{V}\Big)_\nu\Big)^2\,d\sigma\Big/\int_{\partial M}V\Big(\frac{z}{V}\Big)^2\,d\sigma\Big]^{\frac{1}{2}}.
\end{align}
Hence, the desired estimate \eqref{9-th-222} follows from \eqref{9-thm-Proof-7} and \eqref{9-thm-Proof-8}.

\subsection{Proof of Theorem \ref{1-thm10}}
(1) Let $f$ be the first eigenfunction corresponding to the first eigenvalue $p_1$ of the Steklov problem:
$$
\Delta_\phi f-\frac{\Delta_\phi V}{V}f=0\  {\rm on }\ M,\ \ \ V\Big(\frac{f}{V}\Big)_\nu=p \frac{f}{V} \ {\rm on }\ \partial M.
$$
Then from \eqref{2-corr-111}, we have
\begin{align}\label{10-thm-Proof-1}
(m-1)&Kp_1\int_{\partial M}V\Big(\frac{f}{V}\Big)^2\,d\sigma\notag\\
=&(m-1)K\int_{\partial M}V^2\Big(\frac{f}{V}\Big)_\nu\frac{f}{V}\,d\sigma\notag\\
=&(m-1)K\int_MV^2\Big|\nabla\frac{f}{V}\Big|^2\,d\mu\notag\\
\geq&\int_{\partial M}\Big[c_1V^2\Big|\overline{\nabla}\frac{f}{V}\Big|^2+c_2V^3\Big(\Big(\frac{f}{V}\Big)_\nu\Big)^2
+2V^2\Big(\frac{f}{V}\Big)_\nu\Big(\overline{\Delta}_\phi f-\frac{\overline{\Delta}_\phi V}{V}f\Big)\Big]\,d\sigma\notag\\
=&\int_{\partial M}\Big[c_2p_1^2V\Big(\frac{f}{V}\Big)^2+(c_1-2p_1)V^2\Big|\overline{\nabla}\frac{f}{V}\Big|^2\Big]\,d\sigma\notag\\
>&(c_1-2p_1)\int_{\partial M}V^2\Big|\overline{\nabla}\frac{f}{V}\Big|^2\,d\sigma.
\end{align}
Applying the inequality
\begin{align}\label{10-thm-Proof-2}
\eta_1\int_{\partial M}V\Big(\frac{f}{V}\Big)^2\,d\sigma\leq\int_{\partial M}V^2\Big|\overline{\nabla}\frac{f}{V}\Big|^2\,d\sigma,
\end{align}
in \eqref{10-thm-Proof-1} gives
\begin{align}\label{10-thm-Proof-3}
0>[(c_1-2p_1)\eta_1-(m-1)Kp_1]\int_{\partial M}V\Big(\frac{f}{V}\Big)^2\,d\sigma.
\end{align}
Thus, we obtain \eqref{10-th-111}.

(2) Let $f$ be a solution of the following equation:
$$\Delta_\phi f-\frac{\Delta_\phi V}{V}f=0\  {\rm on }\ M,\ \ \ -\beta\Big(\overline{\Delta}_\phi f-\frac{\overline{\Delta}_\phi V}{V} f\Big)+V\Big(\frac{f}{V}\Big)_\nu=\lambda \frac{f}{V} \ {\rm on }\ \partial M,$$
then on $\partial M$, we have
$$\overline{\Delta}_\phi f-\frac{\overline{\Delta}_\phi V}{V} f=\frac{1}{\beta}V\Big(\frac{f}{V}\Big)_\nu-\frac{\lambda_{1,\beta}}{\beta} \frac{f}{V}.$$
Then \eqref{2-corr-111} becomes
\begin{align}\label{10-thm-Proof-4}
0\geq&\int_{\partial M}\Big[c_1V^2\Big|\overline{\nabla}\frac{f}{V}\Big|^2+c_2V^3\Big(\Big(\frac{f}{V}\Big)_\nu\Big)^2
+2V^2\Big(\frac{f}{V}\Big)_\nu\Big(\overline{\Delta}_\phi f-\frac{\overline{\Delta}_\phi V}{V}f\Big)\Big]\,d\sigma\notag\\
=&\int_{\partial M}\Big[\Big(c_2+\frac{2}{\beta}\Big)V^3\Big(\Big(\frac{f}{V}\Big)_\nu\Big)^2-\frac{c_1
+2\lambda_{1,\beta}}{\beta}V^2\Big(\frac{f}{V}\Big)_\nu \frac{f}{V}+\frac{c_1\lambda_{1,\beta}}{\beta}V\Big(\frac{f}{V}\Big)^2
\Big]\,d\sigma.
\end{align}
Inserting the inequality
$$\Big(c_2+\frac{2}{\beta}\Big) V^3\Big(\Big(\frac{f}{V}\Big)_\nu\Big)^2
-\frac{c_1+2\lambda_{1,\beta}}{\beta}V^2\Big(\frac{f}{V}\Big)_\nu \frac{f}{V}
\geq-\frac{(c_1+2\lambda_{1,\beta})^2}{4\beta(2+\beta c_2)}
V\Big(\frac{f}{V}\Big)^2$$
into \eqref{10-thm-Proof-4} yields
\begin{align}\label{10-thm-Proof-5}
0\geq\Big[\frac{c_1\lambda_{1,\beta}}{\beta}-\frac{(c_1+2\lambda_{1,\beta})^2}{4\beta(2+\beta c_2)}\Big]\int_{\partial M}V\Big(\frac{f}{V}\Big)^2\,d\sigma,
\end{align}
which shows that
\begin{align}\label{10-thm-Proof-6}
4(2+\beta c_2)c_1\lambda_{1,\beta}-(c_1+2\lambda_{1,\beta})^2\leq0.
\end{align}
Solving this quadratic inequality with respect to $\lambda_{1,\beta}$, we have that either
\begin{align}\label{10-thm-Proof-7}
\lambda_{1,\beta}\geq\frac{c_1}{2}\Big(1+\beta c_2+\sqrt{\beta^2 c_2^2+2\beta c_2}\Big),
\end{align}
or
\begin{align}\label{10-thm-Proof-8}
\lambda_{1,\beta}\leq\frac{c_1}{2}\Big(1+\beta c_2-\sqrt{\beta^2 c_2^2+2\beta c_2}\Big).
\end{align}

Note that
\begin{align}\label{10-thm-Proof-9}
0=&-\int_MV\frac{f}{V}\Big(\Delta_\phi f-\frac{\Delta_\phi V}{V}f\Big) \,d\mu\notag\\
=&\int_MV^2\Big|\nabla\frac{f}{V}\Big|^2\,d\mu-\int_{\partial M}V^2\Big(\frac{f}{V}\Big)_\nu\frac{f}{V}\,d\sigma\notag\\
=&\int_MV^2\Big|\nabla\frac{f}{V}\Big|^2\,d\mu-\lambda_{1,\beta}\int_{\partial M}V\Big(\frac{f}{V}\Big)^2\,d\sigma+\beta\int_{\partial M}V^2\Big|\overline{\nabla}\frac{f}{V}\Big|^2\,d\sigma.
\end{align}
By the Rayleigh-Ritz formula(or see \eqref{9-thm-Proof-8} with $\beta=0$), we have
\begin{align}\label{10-thm-Proof-10}
\int_MV^2\Big|\nabla\frac{f}{V}\Big|^2\,d\mu\geq p_1\int_{\partial M}V\Big(\frac{f}{V}\Big)^2\,d\sigma
\end{align}
and
\begin{align}\label{10-thm-Proof-11}
\int_{\partial M}V^2\Big|\overline{\nabla}\frac{f}{V}\Big|^2\,d\sigma\geq\eta_1\int_{\partial M}V\Big(\frac{f}{V}\Big)^2\,d\sigma.
\end{align}
Inserting \eqref{10-thm-Proof-10} and \eqref{10-thm-Proof-11} into \eqref{10-thm-Proof-9} gives
\begin{align}\label{10-thm-Proof-12}
\lambda_{1,\beta}\geq p_1+\beta\eta_1>\frac{c_1}{2}+\beta c_1c_2,
\end{align}
which shows that \eqref{10-thm-Proof-8} does not occur.
Therefore, we have that \eqref{10-th-222} holds and the proof of
Theorem \ref{1-thm10} is completed.

\bibliographystyle{Plain}

\end{document}